\documentclass[12pt]{amsart}
\newtheorem{theorem}{Theorem}[section]
\newtheorem{proposition}[theorem]{Proposition}
\newtheorem{corollary}[theorem]{Corollary}
\newtheorem{definition}[theorem]{Definition}

\begin{document}

\baselineskip=16pt

\title[Reduction of stable principal bundles]{Stable principal
bundles and reduction of structure group}

\author[I. Biswas]{Indranil Biswas}

\address{School of Mathematics, Tata Institute of Fundamental
Research, Homi Bhabha Road, Mumbai 400005, India}

\email{indranil@math.tifr.res.in}

\date{}

\begin{abstract}

Let $E_G$ be a stable principal $G$--bundle over a compact
connected K\"ahler manifold, where $G$ is a connected
reductive linear algebraic group defined over $\mathbb C$.
Let $H\subset G$ be a complex reductive subgroup which is
not necessarily connected, and let $E_H\subset E_G$
be a holomorphic reduction of structure group. We
prove that $E_H$ is preserved by the Einstein--Hermitian
connection on $E_G$. Using this we show that if $E_H$ is
a minimal reductive reduction in the sense that there
is no complex reductive proper subgroup of $H$ to which
$E_H$ admits a holomorphic reduction of structure group,
then $E_H$ is unique in the following sense: For any
other minimal reductive reduction $(H'\, , E_{H'})$ of
$E_G$, there is some element $g\in G$ such that $H'=
g^{-1}Hg$ and $E_{H'}= E_Hg$. As an application, we give
an affirmative answer to a question posed in \cite{BK}.

\end{abstract}

\maketitle

\section{Introduction}

In \cite{BK}, the notion of algebraic holonomy for a stable
vector bundle defined over a normal complex projective variety
is introduced. Let $X$ be a simply connected smooth
complex projective
variety equipped with a K\"ahler--Einstein metric. In
\cite[Question 9]{BK} it is asked whether the 
algebraic holonomy of $TM$ coincides with
the complexification of the differential geometric holonomy
of $M$. The present work started
by trying to answer it. We prove that for any stable
vector bundle $E$ defined over a complex projective manifold,
the algebraic holonomy coincides with the complexification
of the
differential geometric holonomy of the Einstein--Hermitian
connection on $E$. Since a K\"ahler--Einstein connection
is also an Einstein--Hermitian connection, this answers the
above question affirmatively.

Let $M$ be a compact connected K\"ahler manifold equipped
with a K\"ahler form. Let $G$ be a connected reductive
linear algebraic group defined over $\mathbb C$.
It is known that any stable principal $G$--bundle
over $M$ admits a unique Einstein--Hermitian connection.

Our main result is the following (Theorem \ref{thm1}):

\begin{theorem}\label{th.0}
Let $E_G$ be a stable $G$--bundle over $M$. Let 
$H$ be a complex reductive subgroup of $G$ which is not
necessarily connected, and let $E_H\, \subset\, E_G$ be a
holomorphic reduction of structure group of $E_G$ to $H$.
Then the Einstein--Hermitian connection on $E_G$ is
induced by a connection on $E_H$.
\end{theorem}

Let $H$ be a complex reductive subgroup of $G$ which is not
necessarily connected, and let $E_H\, \subset\, E_G$ be a
holomorphic reduction of structure group to $H$ of
a stable $G$--bundle $E_G$ defined over $M$. Assume that
there is no complex reductive proper subgroup of $H$ 
to which $E_H$ admits a holomorphic reduction of structure
group. Such reductions will be called the minimal reductive
ones. Theorem \ref{th.0} says that $E_H$ is preserved
by the Einstein--Hermitian connection on $E_G$.

Fix a point $x_0\, \in\, M$, and also choose a point
in the fiber
$z\, \in\, (E_G)_z$. Taking parallel translations, for the
Einstein--Hermitian connection on $E_G$, of $z$
along piecewise smooth paths in $M$ based at $x_0$ we get a
subset of $E_G$. The topological closure, in $E_G$, of this subset
gives a smooth reduction of structure group of $E_G$ to a
compact subgroup $\overline{K}_{E}\, \subset\, G$. The
corresponding smooth reduction of structure group
$E^{z_0}_{\overline{K}^{\mathbb C}_{E}}$ of $E_G$ to the
Zariski closure $\overline{K}^{\mathbb C}_{E}$ of
$\overline{K}_{E}$ in $G$ is actually holomorphic.
We prove the following (see Theorem \ref{thm2}):

\begin{theorem}\label{th.1}
There is a point $z_0\, \in\, (E_G)_z$ such that the
above minimal reductive reduction $E_H$ coincides with
$E^{z_0}_{\overline{K}^{\mathbb C}_{E}}$. In particular,
$H$ coincides with $\overline{K}^{\mathbb C}_{E}$ for
such a base point $z_0$.
\end{theorem}

Consequently, if $(H\, , E_H)$ and $(H'\, , E_{H'})$
are two minimal reductive reductions of $E_G$, then there
is an element $g\, \in\, G$ such that $H'\, =\, g^{-1}Hg$
and $E_{H'}\, =\, E_Hg$. 
This gives us the following corollary (Corollary
\ref{coro2}):

\begin{corollary}\label{co.1}
Let $E_G$ be a stable principal $G$--bundle over $M$. Then
there is a unique holomorphic sub--fiber bundle
${\mathcal G}_{E_G}$ of the adjoint bundle
${\rm Ad}(E_G)$, with fibers being subgroups, that satisfies the
following condition: For any minimal reductive
reduction $E_H$ of $E_G$, the adjoint bundle ${\rm Ad}(E_H)$,
which is a sub--fiber bundle of ${\rm Ad}(E_G)$, 
coincides with ${\mathcal G}_{E_G}$.
\end{corollary}

The sub--fiber
bundle ${\mathcal G}_{E_G}\, \subset\, {\rm Ad}(E_G)$
in Corollary \ref{co.1} is the
complexification of the holonomy of the Einstein--Hermitian
connection on $E_G$. While the stability condition does
not depend on the choice of a K\"ahler form in a given K\"ahler
class, the Einstein--Hermitian connection on a stable bundle
depends on the choice of the K\"ahler form. We note that the
sub--fiber bundle ${\mathcal G}_{E_G}$ does not depend either
on the K\"ahler form or on the K\"ahler class as long as
$E_G$ remains stable. When $M$ is a complex projective manifold,
${\mathcal G}_{E_G}$ coincides with the algebraic monodromy
of $E_G$ introduced in \cite{BK}.

\section{Einstein--Hermitian connection and reduction of
structure group}

Let $M$ be a compact connected K\"ahler manifold equipped with a
K\"ahler form $\omega$. The \textit{degree} of a torsionfree
coherent analytic sheaf $V$ defined over a dense open subset
$U\, \subset\, M$, whose complement $M\setminus U$ is an analytic
subset of $M$ of complex codimension at least two,
is defined to be
\[
\text{degree}(V) \, :=\, \int_M c_1(\iota_* V)
\omega^{\dim_{\mathbb C} M-1} \, \in\, {\mathbb R}\, ,
\]
where $\iota \, :\, U \, \hookrightarrow\, M$ is the inclusion
map; the codimension condition ensures that $\iota_* V$ is a
coherent analytic sheaf.

Let $G$ be a connected reductive linear
algebraic group defined over $\mathbb C$. Since the degree has
been defined, we have the notions of a stable
$G$--bundle and a polystable $G$--bundle over $M$.
See \cite{RS},
\cite{AB}, \cite{Ra} for definitions of stable and polystable 
principal $G$--bundles.

Fix a maximal torus $T\, \subset\, G$ and a 
Borel subgroup $B\, \subset\, G$ containing $T$. We also
fix a maximal compact subgroup $K$ of $G$. By a parabolic
subgroup of $G$ we will mean one containing $B$. Therefore,
the Levi quotient $L(Q)$ of any parabolic subgroup $Q\,\subset\,
G$ is also a subgroup of $Q$. Let $Z_0(G)$ denote the
connected component of the center of $G$ that contains
the identity element.

Take any polystable principal $G$--bundle $E_G$ over $M$.
Therefore, there is a Levi subgroup $L(P)$ associated to some
parabolic subgroup $P\, \subset\, G$ (the subgroup $P$ need
not be proper) along with a holomorphic reduction of structure
group $E_{L(P)}\, \subset\, E_G$ to $L(P)$, such that
\begin{itemize}
\item the principal $L(P)$--bundle $E_{L(P)}$ is stable, and

\item the principal
$P$--bundle $E_P\, :=\, E_{L(P)}(P)$, obtained
by extending the structure group of
$E_{L(P)}$ using the inclusion $L(P)$ in $P$, is an admissible
reduction of $E_G$.
\end{itemize}
Note that since $L(P)\, \subset\,P \, \subset\,G$,
the $P$--bundle $E_P$ is a reduction of structure group of $E_G$
to $P$. The condition that $E_P$ is an admissible reduction
of structure group of $E_G$ means that for each character $\chi$
of $P$, which is trivial on $Z_0(G)$, the
associated line bundle $E_P(\chi)$ over $M$ is of degree zero.
The $G$--bundle $E_G$ is stable if and only if $P\, =\, G$.

Since $E_G$ and $E_{L(P)}$ are polystable, they admit
unique Einstein--Hermitian connections \cite[p. 208,
Theorem 0.1]{AB}. It should be clarified that the
Einstein--Hermitian reduction of structure to a maximal compact
subgroup depends on the choice of the maximal compact
subgroup. Even for a fixed maximal compact
subgroup, the Einstein--Hermitian reduction of structure
group need not be unique.
However, once the K\"ahler form on $M$ is fixed,
the Einstein--Hermitian connection on
a polystable principal bundle over $M$ is unique. Let
$\nabla^{L(P)}$ be the Einstein--Hermitian connection on the
stable $L(P)$--bundle $E_{L(P)}$. Let $\nabla^G$ be the
connection on $E_G$ induced by $\nabla^{L(P)}$.
Note that a connection on $L(P)$
induces a connection on any fiber bundle associated to $E_{L(P)}$.
The principal $G$--bundle $E_G$ is associated to
$E_{L(P)}$ for the left translation action of $L(P)$ on $G$, and
$\nabla^G$ is the connection on it obtained from $\nabla^{L(P)}$.

\begin{proposition}\label{prop.1}
The induced connection $\nabla^G$ on $E_G$ coincides with the
unique Einstein--Hermitian connection on $E_G$.
\end{proposition}

\begin{proof}
Since the inclusion map $L(P)\, \hookrightarrow\, G$ need not
take the connected component, containing the identity
element, of the center of $L(P)$ into $Z_0(G)$, the proposition
does not follow immediately from \cite[p. 208, Theorem 0.1]{AB}.
The connection $\nabla^{L(P)}$ on $E_{L(P)}$ is induced
by a connection on a smooth reduction of structure group of
$E_{L(P)}$ to a maximal compact subgroup of $L(P)$ (this is a part
of the definition of a Einstein--Hermitian connection). A maximal
compact subgroup of $L(P)$ is contained in a maximal compact
subgroup of $G$. Consequently, the connection $\nabla^G$ on $E_G$
is induced
by a connection on a smooth reduction of structure group of
$E_G$ to a maximal compact subgroup of $G$. Therefore, to prove
that $\nabla^G$ is the Einstein--Hermitian connection on $E_G$ it
suffices to show that $\nabla^G$ satisfies the Einstein--Hermitian
equation.

Let ${\mathfrak z}({\mathfrak l}(P))$ be the Lie algebra
of the center of $L(P)$. For any $\theta\, \in\,
{\mathfrak z}({\mathfrak l}(P))$, the holomorphic
section of the adjoint bundle $\text{ad}(E_{L(P)})$ given by
$\theta$ will be denoted by $\widehat{\theta}$; the vector
bundle $\text{ad}(E_{L(P)})$ is associated to $E_{L(P)}$ for
the adjoint action of $L(P)$ on its Lie algebra.
Let $\Lambda_\omega$ be the adjoint of multiplication
by $\omega$ of differential form on $M$.
That the connection
$\nabla^{L(P)}$ is Einstein--Hermitian means that there is
an element $\theta \, \in\, {\mathfrak z}({\mathfrak l}(P))$
such that the Einstein--Hermitian equation
\begin{equation}\label{eq1}
\Lambda_\omega {\mathcal K}_{\nabla^{L(P)}} \, =\,
\widehat{\theta}
\end{equation}
holds.

If $\theta$ in eqn. \eqref{eq1} is in the Lie algebra
${\mathfrak z}(G)$ of $Z_0(G)$, then $\nabla^G$ is a
Einstein--Hermitian connection on $E_G$. Therefore, in that
case the proposition is proved. Assume that
\begin{equation}\label{notin}
\theta\, \not\in\, {\mathfrak z}(G)\, .
\end{equation}
Fix a character $\chi$ of $P$ which is trivial on $Z_0(G)$ but
satisfies the following condition: the homomorphism of Lie
algebras
\begin{equation}\label{eqb1}
{\rm d}\chi\, :\, {\mathfrak p}\, \longrightarrow\, {\mathbb C}
\end{equation}
given by $\chi$, where $\mathfrak p$ is the Lie algebra of
$P$, is nonzero on $\theta$. (It is easy to check that the
group of characters of $P$ coincides with the group of
characters of $L(P)$.)

Consider the holomorphic line bundle $L_\chi \, :=\,
E_P(\chi)$ over $M$
associated to $E_P$ for $\chi$. Let $\nabla^\chi$
be the connection on $L_\chi$ induced by the connection
on $E_P$ given by $\nabla^{L(P)}$. Since $\nabla^{L(P)}$
is a Einstein--Hermitian connection, the connection
$\nabla^\chi$ is also Einstein--Hermitian. Indeed, if
${\mathcal K}(\nabla^\chi)$ is the curvature of
$\nabla^\chi$, then
\[
\Lambda_\omega {\mathcal K}(\nabla^\chi)\, =\,
{\rm d}\chi (\theta)\, ,
\]
where ${\rm d}\chi$ is the homomorphism in eqn. \eqref{eqb1}
and $\theta$ is the element in eqn. \eqref{eq1}. Therefore,
\begin{equation}\label{eq2}
\text{degree}(L_\chi) \, =\, \frac{{\rm d}\chi
(\theta)\sqrt{-1}}{2{\pi}d} \int_M \omega^d\, ,
\end{equation}
where $d \, =\, \dim_{\mathbb C}M$; see
\cite[p. 103, Proposition 2.1]{Ko}.

The condition
that $E_P\, \subset\, E_G$ is an admissible reduction of
structure group says that
\[
\text{degree}(L_\chi) \, =\, 0\, .
\]
This, in view of the assumption in eqn. \eqref{notin}
that ${\rm d}\chi (\theta)\,
\not=\, 0$, contradicts eqn. \eqref{eq2}. Therefore, we conclude
that $\theta\, \in\, {\mathfrak z}(G)$. This immediately implies
that the connection $\nabla^G$ on $E_G$ is
Einstein--Hermitian. This completes the
proof of the proposition.
\end{proof}

Let $E_G$ be a holomorphic principal $G$--bundle over $M$
and $\nabla'$ a $C^\infty$ connection on $E_G$ compatible
with the holomorphic structure of $E_G$. Such a connection
is called a complex connection; see \cite[p. 230, Definition
3.1(1)]{AB} for the precise definition of a complex connection.
A connection $\widetilde \nabla$ is complex if and only if
the $(0\, ,2)$--Hodge type component of the curvature of
$\widetilde \nabla$ vanishes.

\begin{definition}\label{def1}
{\rm Let $H$ be a closed complex subgroup of $G$ and $E_H\,
\subset\, E_G$ a $C^\infty$ reduction of structure group of
$E_G$ to $H$. We will say that $E_H$ is \textit{preserved} by
the connection
$\nabla'$ if $\nabla'$ induces a connection on $E_H$.}
\end{definition}

It follows immediately
that $E_H$ is preserved by $\nabla'$ if and only if there
is a smooth connection on $E_H$ that induces the connection
$\nabla'$ on $E_G$. It is easy to see that $E_H$ is preserved
by $\nabla'$ if and only if for each point $z\, \in\, E_H$, the
horizontal subspace in $T_z E_G$ for the connection
$\nabla'$ is contained in $T_z E_H$.
If $E_H$ is preserved by $\nabla'$, then
$E_H$ is a holomorphic reduction of structure group of $E_G$.

\begin{theorem}\label{thm1}
Let $E_G$ be a stable principal
$G$--bundle over $M$ and $\nabla^G$ the
Einstein--Hermitian connection on $E_G$.
Let $H$ be a complex reductive subgroup of $G$ which is not
necessarily connected, and let $E_H\, \subset\, E_G$
be a holomorphic reduction of structure group of $E_G$
to $H$. Then $E_H$ is preserved by the connection $\nabla^G$.
\end{theorem}

\begin{proof}
Let $H_0\, \subset\, H$ be the connected component containing the
identity element. Set
\[
X \, :=\, E_H/H_0\, ,
\]
which is finite \'etale Galois cover of $M$ with Galois
group $H/H_0$. Let
\begin{equation}\label{def.p}
p\, :\, X\, \longrightarrow\, M
\end{equation}
be the projection.
We note that $p^*E_H$ has a canonical reduction of structure
group to the subgroup $H_0\, \subset\, H$. Set
\[
F_G\, :=\, p^*E_G\, .
\]
Let $F_{H_0}\, \subset\, F_G$ be the reduction of structure group
to $H_0$ obtained from the canonical reduction of structure group
of $p^*E_H$ to $H_0$.

Equip $X$ with the K\"ahler form $p^*\omega$.
The Einstein--Hermitian connection on $E_G$ pulls back to a
Einstein--Hermitian connection on $F_G$ for the K\"ahler form
$p^*\omega$. Therefore, $F_G$ is polystable with respect to
$p^*\omega$. Let $\text{ad}(F_G)$ be the adjoint bundle over $X$.
We recall that $\text{ad}(F_G)$ is associated to $F_G$ for the
adjoint action of $G$ on its Lie algebra. The connection on
the vector bundle
$\text{ad}(F_G)$ induced by the Einstein--Hermitian connection
of $F_G$ is clearly Einstein--Hermitian. Hence the adjoint
vector bundle $\text{ad}(F_G)$ is polystable.

Let $Z(G)$ be the center of $G$. Set 
\[
H'\, :=\, H_0/(H_0\cap Z(G))\, .
\]
Therefore, $H'$ is a complex reductive subgroup of the complex
semisimple group $G'\, :=\, G/Z(G)$.
Let ${\mathfrak g}$ (respectively,
${\mathfrak h}_0$) be the Lie algebra of $G$ (respectively, $H_0$).
Note that the adjoint action makes ${\mathfrak g}$ (respectively,
${\mathfrak h}_0$) a $G'$--module (respectively, $H'$--module).
We will also consider ${\mathfrak g}$ as a $H'$--module using
the inclusion of $H'$ in $G'$.

Since ${\mathfrak g}$ is a faithful $G'$--module, and $H'$ is
a reductive
subgroup of $G'$, there is a positive integer $N$ and nonnegative
integers $a_i, b_i$, $i\, \in\, [1\, ,N]$, such
that the $H'$--module
${\mathfrak h}_0$ is a direct summand of the $H'$--module
\begin{equation}\label{vec.s.}
\bigoplus_{i=1}^N {\mathfrak g}^{\otimes a_i}\otimes
({\mathfrak g}^*)^{\otimes b_i}\, =\,
\bigoplus_{i=1}^N {\mathfrak g}^{\otimes a_i}\otimes
{\mathfrak g}^{\otimes b_i}
\end{equation}
\cite[p. 40, Proposition 3.1]{De}; since $H'$ is complex
reductive, any exact sequence of $H'$--modules splits; also,
${\mathfrak g}\, = \, {\mathfrak g}^*$ as $G$ is reductive.
Therefore, the adjoint vector bundle $\text{ad}(F_{H_0})$
is a direct summand of the vector bundle
\begin{equation}\label{vec.bu.s.}
\bigoplus_{i=1}^N \text{ad}(F_G)^{\otimes a_i}\otimes
\text{ad}(F_G)^{\otimes b_i}\, .
\end{equation}
We note that the vector bundle in eqn. \eqref{vec.bu.s.}
is associated to $F_{H_0}$ for the $H_0$--module in
eqn. \eqref{vec.s.}.

Since the adjoint vector bundle $\text{ad}(F_G)$ is polystable
of degree zero (recall that $\text{ad}(F_G)\, =\,\text{ad}(F_G)^*$),
the vector bundle
\[
\text{ad}(F_G)^{\otimes a}\otimes
\text{ad}(F_G)^{\otimes b}
\]
is polystable of degree zero for all $a, b\, \in\, {\mathbb N}$
\cite[p. 224, Theorem 3.9]{AB}.
Therefore, the vector bundle in eqn. \eqref{vec.bu.s.} is
polystable of degree zero. Since $\text{ad}(F_{H_0})$ is a direct
summand of it of degree zero, we conclude that
$\text{ad}(F_{H_0})$ is also polystable. Consequently, the
principal $H_0$--bundle $F_{H_0}$ is polystable
\cite[p. 224, Corollary 3.8]{AB}. Hence $F_{H_0}$ admits a unique
Einstein--Hermitian connection \cite[p. 208,
Theorem 0.1]{AB}. Let $\nabla^{H_0}$ be the Einstein--Hermitian
connection on $F_{H_0}$. So, there is an element $\nu \, \in\,
{\mathfrak h}$ such that
\begin{equation}\label{ehc}
\Lambda_{p^*\omega} {\mathcal K}(\nabla^{H_0}) \, =\,
\widehat{\nu}\, ,
\end{equation}
where ${\mathcal K}(\nabla^{H_0}) $ is the curvature of
$\nabla^{H_0}$, and $\widehat{\nu}$ is the holomorphic section
of $\text{ad}(p^*E_H)\, =\, \text{ad}(F_{H_0})$ given by $\nu$.

Since $H/H_0$ is a finite group, giving a connection on a
principal $H_0$--bundle is equivalent to giving a connection on
the principal $H$--bundle obtained from it
by extension of structure
group. From the uniqueness of the Einstein--Hermitian connection
on $F_{H_0}$ it follows that the corresponding connection on
$p^*E_H$ is left invariant by the action of the Galois group
$H/H_0$.
Consequently, the connection on $p^*E_H$ given by $\nabla^{H_0}$
descends to a connection on $E_H$. Let $\nabla^H$ denote the
connection on $E_H$ obtained this way. It is clear that
$\nabla^H$ is a complex connection.

Let $\nabla$ be the complex connection on $E_G$ induced by the
above connection $\nabla^H$ on $E_H$. We will first show that
$\nabla$ is unitary, which means that
$\nabla$ is induced by connection
on a smooth reduction of structure group of
$E_G$ to a maximal compact subgroup of $G$. Then we will show
that $\nabla$ satisfies the Einstein--Hermitian equation.

To prove that $\nabla$ is unitary, fix a point $x_0\, \in\, M$,
and also fix a point $z_0\, \in\, (E_G)_{x_0}$ in the fiber
of $E_G$ over $x_0$. Taking parallel translations of $z_0$,
with respect to the connection $\nabla$, along
piecewise smooth
paths in $M$ based at $x_0$ we get a subset ${\mathcal S}$
of $E_G$. Sending any $g\, \in\, G$ to the point
$z_0g\, \in\, (E_G)_{x_0}$ we get an isomorphism
$G\, \longrightarrow\, (E_G)_{x_0}$. Using this isomorphism,
the intersection $(E_G)_{x_0}\bigcap {\mathcal S}$ is a 
subgroup $K_E$ of $G$. The condition that the connection $\nabla$
is unitary is equivalent to the condition that $K_E$ is contained
in some compact subgroup of $G$.

Fix a point $x\,\in\, p^{-1}(x_0)$, and also fix a point $z\, \in
\, (p^*E_G)_x$ that projects to $z$, where $p$ is the covering map
in eqn. \eqref{def.p}. Consider parallel translations of $z$,
with respect to the connection $p^*\nabla$ on $p^*E_G \, =:\, F_G$,
along piecewise smooth paths in $X$ based at $x$. As before, we
get a subgroup $K_F\, \subset\, G$ from the resulting
subset of $F_G$. It is easy to see that $K_F$ is a finite index
subgroup of the group $K_E$ constructed above.

We note that the connection $p^*\nabla$ on $F_G$ 
is induced by the unitary connection $\nabla^{H_0}$ on
the reduction $F_{H_0}$ of $F_G$. Therefore, the
connection $p^*\nabla$ is unitary. Consequently,
the subgroup $K_F\, \subset\, G$ is contained in a compact
subgroup of $G$. Using this together with the
observation that $K_F$ is a finite index
subgroup of the subgroup $K_E\, \subset\, G$ we conclude that
$K_E$ is also contained in a compact subgroup of $G$.
Thus the connection $\nabla$ is unitary.

We will now show that $\nabla$ satisfies the
Einstein--Hermitian equation.

Let ${\mathcal K}(\nabla)$
be the curvature of $\nabla$. From eqn. \eqref{ehc}
it follows immediately that $\Lambda_{p^*\omega}
{\mathcal K}(\nabla^{H_0})$ is a holomorphic section of
$\text{ad}(p^*E_H)$. Consequently, 
$\Lambda_{\omega} {\mathcal K}(\nabla)$ is a holomorphic
section of the adjoint vector bundle $\text{ad}(E_G)$.
Since $E_G$ is stable, all holomorphic sections of
$\text{ad}(E_G)$ are given by the Lie algebra
${\mathfrak z}(G)$ of $Z_0(G)$ \cite[Proposition 3.3]{BG}.
In other words, there
is an element $\theta\, \in\, {\mathfrak z}(G)$
such that $\Lambda_{\omega} {\mathcal K}(\nabla)$
coincides with the
section of $\text{ad}(E_G)$ given by $\theta$.
Thus, we conclude that the connection
$\nabla$ is the unique Einstein--Hermitian connection on $E_G$.
Since $\nabla$ is induced by a connection on $E_H$, the proof
of the theorem is complete.
\end{proof}

Proposition \ref{prop.1} and Theorem \ref{thm1} together have
the following corollary.

\begin{corollary}\label{cor1}
Let $E_G$ be a polystable $G$--bundle over $M$. Take
$E_{L(P)}$ as in Proposition \ref{prop.1}. Let
$H$ be a complex reductive subgroup of $L(P)$ (not
necessarily connected), and let $E_H\, \subset\, E_{L(P)}$
be a holomorphic reduction of structure group
of $E_{L(P)}$ to $H$. Then $E_H$ is preserved by the
Einstein--Hermitian connection on $E_G$.
\end{corollary}

\section{Properties of a smallest reduction}

Let $E_G$ be a stable principal $G$--bundle over $M$.
Fix a point $x_0\, \in\, M$, and also fix a point
$z_0\, \in\, (E_G)_{x_0}$ in the fiber over $x_0$. Let
$\nabla^G$ be the Einstein--Hermitian connection on $E_G$.

Taking parallel translations of $z_0$, with respect
to $\nabla^G$, along piecewise smooth
paths in $M$ based at $x_0$ we get a subset ${\mathcal S}$
of $E_G$. Let $\overline{\mathcal S}$ be the topological
closure of $\mathcal S$ in $E_G$.
The map $G\, \longrightarrow\, (E_G)_{x_0}$
defined by $g\, \longmapsto\, z_0g$ is an 
isomorphism. Using this isomorphism, the intersection
$(E_G)_{x_0}\bigcap \overline{\mathcal S}$ gives a compact
subgroup $\overline{K}_{E} \, \subset\, G$. The subset
$\overline{\mathcal S}\, \subset\, E_G$ is a $C^\infty$
reduction of structure group of $E_G$ to the subgroup
$\overline{K}_{E}$. Let
$\overline{K}^{\mathbb C}_{E}$ be the complex reductive
subgroup of $G$ obtained by taking the Zariski closure
of $\overline{K}_{E}$ in $G$. The subset
\begin{equation}\label{mon.red.}
E^{z_0}_{\overline{K}^{\mathbb C}_{E}}\, :=\,
{\mathcal S}\overline{K}^{\mathbb C}_{E}\, \subset\, E_G
\end{equation}
is a holomorphic reduction of structure group of $E_G$
to $\overline{K}^{\mathbb C}_{E}$; see \cite[Section 3]{Bi}
for the details. It is easy to see that the principal
$\overline{K}^{\mathbb C}_{E}$--bundle
$E^{z_0}_{\overline{K}^{\mathbb C}_{E}}$ is the extension
of structure group of the principal $\overline{K}_{E}$--bundle
$\overline{\mathcal S}$.

We note that in \cite{Bi}, the point $z_0$ is taken to be
in the subset of $E_G$ given by an Einstein--Hermitian
reduction of structure group (see \cite[p. 71, (3.20)]{Bi}).
This was only to make $\overline{K}_{E}$ lie inside a fixed
maximal compact subgroup of $G$.
If we replace the base point $z_0$ by $z_0g$, where $g$ is any
point of $G$,
then it is easy to see that the subset $\mathcal S$
constructed above using $z_0$ gets replaced by ${\mathcal S}g$.
Therefore, the subgroup $\overline{K}^{\mathbb C}_{E}$
in eqn. \eqref{mon.red.} gets replaced by 
$g^{-1}\overline{K}^{\mathbb C}_{E}g$, and
the reduction $E^{z_0}_{\overline{K}^{\mathbb C}_{E}}
\,\subset\, E_G$ in eqn. \eqref{mon.red.} gets replaced by
$E^{z_0}_{\overline{K}^{\mathbb C}_{E}}g$.

We note that $(\overline{K}^{\mathbb C}_{E}\, ,
E^{z_0}_{\overline{K}^{\mathbb C}_{E}})$ in eqn.
\eqref{mon.red.} is a minimal complex reduction of $E_G$
in the following sense: There is no complex proper
subgroup of $\overline{K}^{\mathbb C}_{E}$ to which
$E^{z_0}_{\overline{K}^{\mathbb C}_{E}}$ admits a holomorphic
reduction of structure group which is preserved by the
connection $\nabla^G$. Indeed, this follows
immediately from the construction of
$E^{z_0}_{\overline{K}^{\mathbb C}_{E}}$. Furthermore,
if $E_{H'}\, \subset\, E_G$ is a holomorphic reduction of
structure group of $E_G$, to a complex subgroup
$H'\, \subset\, G$, satisfying the two conditions:
\begin{itemize}
\item $E_{H'}$ is preserved by $\nabla^G$, and

\item there is no complex proper subgroup of
$H'$ to which $E_{H'}$ admits a holomorphic reduction
of structure group which is also preserved by $\nabla^G$,
\end{itemize}
then it follows immediately that there is a point
$z_0\, \in\, (E_G)_{x_0}$ such that
$E_{H'}\, =\, E^{z_0}_{\overline{K}^{\mathbb C}_{E}}$.
Indeed, $z_0$ can be taken to be any point of
the fiber $(E_{H'})_{x_0}$.

Let $H$ be a complex reductive subgroup of $G$ which is not
necessarily connected, and let $E_H\, \subset\, E_G$ be a
holomorphic reduction of structure group to $H$ of
the stable $G$--bundle $E_G$. This reduction $E_H$ will 
be called a \textit{minimal reductive reduction} of $E_G$
if there is no complex reductive proper subgroup of $H$ 
to which $E_H$ admits a holomorphic reduction of structure
group.

In view of the above observations, using Theorem \ref{thm1}
get the following theorem:

\begin{theorem}\label{thm2}
Let $E_G$ be a stable principal $G$--bundle over a
compact connected K\"ahler manifold $M$, where $G$ is a
connected reductive complex linear algebraic group. Let
$E_H\, \subset\,E_G$ be a minimal reductive
reduction of $E_G$. Fix a point $x_0\, \in\, M$.
Then there is a point $z_0$ in the
fiber $(E_G)_{x_0}$ such that
$E^{z_0}_{\overline{K}^{\mathbb C}_{E}}$ (defined
in eqn. \eqref{mon.red.}) coincides with $E_H$.
In particular, the subgroup
$H$ coincides with $\overline{K}^{\mathbb C}_{E}$ for
such a base point $z_0$.
\end{theorem}

Let $\text{Ad}(E_G)$ be the adjoint bundle for $E_G$.
So $\text{Ad}(E_G)$ is the fiber bundle over $M$ associated to
$E_G$ for the adjoint action of $G$ on itself. The fibers
of $\text{Ad}(E_G)$ are groups isomorphic to $G$.
Let $\text{Ad}(E^{z_0}_{\overline{K}^{\mathbb C}_{E}})$
be the adjoint bundle of the principal
$\overline{K}^{\mathbb C}_{E}$--bundle
$E^{z_0}_{\overline{K}^{\mathbb C}_{E}}$ defined
in eqn. \eqref{mon.red.}. Since $E^{z_0}_{\overline{K}^{\mathbb
C}_{E}}$ is a holomorphic reduction of structure group of
$E_G$, the fiber bundle $\text{Ad}(E^{z_0}_{\overline{K}^{\mathbb
C}_{E}})$ is a holomorphic
sub--fiber bundle of $\text{Ad}(E_G)$ with the fibers
of $\text{Ad}(E^{z_0}_{\overline{K}^{\mathbb
C}_{E}})$ being subgroups of the fibers of $\text{Ad}(E_G)$.

We noted earlier that if $z_0$ is replaced by $z_0g$, where
$g\, \in\, G$, then
the subgroup $\overline{K}^{\mathbb C}_{E}$ gets replaced by
$g^{-1}\overline{K}^{\mathbb C}_{E}g$, and
$E^{z_0}_{\overline{K}^{\mathbb C}_{E}}$ gets replaced
by $E^{z_0}_{\overline{K}^{\mathbb C}_{E}}g$. From this
it follows immediately that the subbundle
\[
\text{Ad}(E^{z_0}_{\overline{K}^{\mathbb
C}_{E}})\, \subset\, \text{Ad}(E_G)
\]
is independent of the choice of the base point $z_0$.

Therefore, from Theorem \ref{thm2}
we have the following corollary:

\begin{corollary}\label{coro2}
Let $E_G$ be a stable principal $G$--bundle over a
compact connected K\"ahler manifold $M$, where $G$ is a
connected reductive complex linear algebraic group. Then there
is a unique holomorphic sub--fiber bundle ${\mathcal G}_{E_G}$
of the adjoint bundle ${\rm Ad}(E_G)$, with
fibers being subgroups, that satisfies the
following condition: For any minimal reductive reduction
$E_H$ of $E_G$, the adjoint bundle ${\rm Ad}(E_H)$,
which is a sub--fiber bundle of ${\rm Ad}(E_G)$, 
coincides with ${\mathcal G}_{E_G}$.
\end{corollary}

{}From Theorem \ref{thm2} it follows that the sub--fiber
bundle ${\mathcal G}_{E_G}\, \subset\, {\rm Ad}(E_G)$
in Corollary \ref{coro2}
is the complexification of the holonomy of the Einstein--Hermitian
connection on $E_G$. We note that for defining stable bundles we
need only the class $[\omega]\, \in\, H^2(M, \, {\mathbb R})$.
In other words, the stability condition does not depend on the
choice of the K\"ahler form in a given K\"ahler class. On the
other hand, the Einstein--Hermitian connection on a stable bundle
depends on the choice of the K\"ahler form. From Corollary
\ref{coro2} it follows that the sub--fiber
bundle ${\mathcal G}_{E_G}$ does not depend on the choice
of the K\"ahler form in a given K\"ahler class. In fact,
it does not even depend on the K\"ahler class as long
as $E_G$ remains stable.

Using Corollary \ref{coro2} and \cite[Theorem 20]{BK}
it follows that the algebraic monodromy of $E_G$
coincides with the
fiber $({\mathcal G}_{E_G})_{x_0}$ when $M$ is a
complex projective manifold. (See \cite{BK}
for algebraic monodromy where it is introduced.)

%%%%%%%%%%%%%%%%%%%%%%%%%%%%%%%%%%%%%%%%%%%%%%%%%%%%%%%%%%%%%%%


\begin{thebibliography}{AAAA}

\bibitem[AB]{AB} B. \textit{Anchouche} and I. \textit{Biswas},
Einstein--Hermitian connections on polystable principal
bundles over a compact K\"ahler manifold, Amer. Jour.
Math. \textbf{123} (2001), 207--228.

\bibitem[BK]{BK} V. \textit{Balaji} and J. \textit{Koll\'ar},
Holonomy groups of stable vector bundles, Preprint,
math.AG/0601120.

\bibitem[BG]{BG} I. \textit{Biswas} and T. L. \textit{G{\'o}mez},
Simplicity of stable principal sheaves, Bull. Lond. Math. Soc.
(to appear).

\bibitem[Bi]{Bi} I. \textit{Biswas}, Stable bundles and
extension of structure group, Diff. Geom. Appl. \textbf{23}
(2005), 67--78.

\bibitem[De]{De} P. \textit{Deligne}, Hodge cycles
on abelian varieties, (in: Hodge cycles, motives,
and Shimura varieties, by P. Deligne,
J. S. Milne, A. Ogus and K.-Y. Shih), Lecture Notes
in Mathematics, 900, Springer-Verlag, Berlin-New York, 1982.

\bibitem[Ko]{Ko} S. \textit{Kobayashi}, Differential Geometry of
Complex Vector Bundles, Publications of the Math. Society of Japan
15, Iwanami Shoten Publishers and Princeton University Press, 1987.

\bibitem[RS]{RS} A. \textit{Ramanathan} and S.
\textit{Subramanian}, Einstein--Hermitian connections
on principal bundles and stability, Jour. Reine
Angew. Math. \textbf{390} (1988), 21--31.

\bibitem[Ra]{Ra} A. \textit{Ramanathan}, Moduli for principal
bundles over algebraic curves: I, Proc. Ind. Acad. Sci.
(Math. Sci.) \textbf{106} (1996), 301--328.

\end{thebibliography}
\end{document}